\documentclass[]{article}
\usepackage{graphicx}
\usepackage{amssymb}
\usepackage{eufrak,mathrsfs,amsmath,longtable}
\newtheorem{theorem}{thm}

\begin{document}
\title{
Geometry of submanifolds of all classes of third-order ODEs as a Riemannian manifold}
\author{Z. Bakhshandeh-Chamazkoti$^\dag$\footnote{Corresponding Author}, A. Behzadi$^{\S}$, R. Bakhshandeh-Chamazkoti$^\ddag$,\\
  and M. Rafie-Rad$^{\S\S}$ \\[3mm]
{\small   $^{\dag, \S, \S\S}$Department of Mathematics, Faculty of Mathematical Sciences, University of Mazandaran,
Babolsar, Iran.}\\[3mm]
{\small $^\ddag$Department of Mathematics, Faculty of Basic Sciences, Babol  Noshirvani University of Technology, Babol, Iran.}\\[3mm]
{\small Emails: $^\dag$z.bakhshandeh.ch@stu.umz.ac.ir,  $^{\S}$behzadi@umz.ac.ir, $^\ddag$r\_bakhshandeh@nit.ac.ir, $^{\S\S}$rafie-rad@umz.ac.ir.}
}
\maketitle
%%%%%%%%%%%%%%%%%%%%%%%
\begin{abstract}
In this paper, we prove that any surface corresponding to  linear second-order ODEs
as a submanifold is minimal in all classes of third-order ODEs $y'''=f(x, y, p, q)$ as a  Riemannian manifold
where $y'=p$ and $y''=q$, if and only if $q_{yy}=0$.
Moreover, we will see the linear second-order ODE with general form $y''=\pm y+\beta(x)$ is the only case that is defined a minimal surface
and is also totally geodesic.
%
%%%%%%%%%%%
\end{abstract}
\vspace{0.7cm}
\section{Introduction}
Riemannian geometry is characterized, and research is oriented towards and shaped
by concepts for examples  geodesics connections, curvature \cite{Jost}.
 Originally, the geometry of submanifolds was only a part of Riemannian geometry
but today it is one of several independent aspects of multi-dimensional generalizations
of the classical theory of surfaces. While Riemannian geometry is the development
of Gauss' idea on intrinsic geometry, the geometry of submanifolds starts
from the idea of the extrinsic geometry of a surface. This theory is devoted to the
study of the position and properties of a submanifold in ambient space, both in their
local and global aspects \cite{Aminov, Anciaux, Willmore}.

On the other hand, the general equivalence problem is about studding when two geometrical objects are mapped
on each other by a certain class of diffeomorphisms. $\rm\acute{E}$lie Cartan developed the general equivalence
problem and provided a systematic procedure for determining the necessary and sufficient conditions.
There are many papers  as applications of Cartan method for third-order ODEs
\cite{Godlinski, Medvedev, Sato} and fourth-order differential operators \cite{Bakhshandeh}.
The theory of moving frames is most closely associated with the name of Cartan, that is
a powerful and algorithmic tool for studying the geometric properties of submanifolds
and their invariants under the action of a transformation group \cite{Olver}.
Now if we apply the moving frame method for Riemannian manifold and geometry of submanifolds, we can obtain very interesting results.
Recently, a lot of research has been carried out about the minimal surfaces in a three-dimensional Riemannian manifold \cite{Daniel, Dillen, Rosenberg}.

The graph of the function $u(x, y)$ is minimal if and only if $u$ satisfies:
\begin{eqnarray}\label{minimal}
(1+u_x^2)\;u_{yy}-2u_xu_yu_{xy}+(1+u_y^2)\;u_{xx}=0.
\end{eqnarray}
The PDE \eqref{minimal} which is non-linear but depends linearly on the second derivatives is called quasi-linear.
We will not discuss the general theory
of existence of solutions to this equation, but rather describe a couple
of special solutions. Of course, a function $u(x, y)$ with vanishing second derivatives
 is solution of Equation \eqref{minimal}. Although these are
not very exciting solutions, the resulting surfaces $\cal S$ being planes, they are
important because of the Calabi theorem \cite{Anciaux}. According to this theorem,
 if $u$ be a solution of Equation \eqref{minimal}, on $U = {\Bbb R}^2$
 such that $\|\nabla u\|_0 < 1$,  (hence $\cal S$ is called
spacelike)  then $u$ is an affine function that means $u$ is a function with vanishing second derivatives.
There is an analogous theorem (which actually predates Calabi's) in the
case of the Riemannian metric $\langle . , .\rangle_0=dx_1^2+dx_2^2+dx_3^2$
and is known as the Bernstein theorem.

It is a classical result that any simply connected minimal surface in Euclidean space ${\Bbb R}^3$ admits a one-parameter family of minimal isometric deformations, called the \emph{associate family}. Conversely, two minimal isometric immersions of the same Riemannian surface into ${\Bbb R}^3$ are associate. These are easy consequences of the Gauss and Codazzi equations in ${\Bbb R}^3$. More generally, analogous results hold for constant mean curvature (CMC) surfaces in $3$-dimensional space forms.

Benoit Daniel in \cite{Daniel} investigated extensions of these results and  related questions for minimal surfaces in the product manifolds
${\Bbb S}^2\times{\Bbb R}$ and  ${\Bbb H}^2\times{\Bbb R}$, where ${\Bbb S}^2$ is the $2$-sphere of curvature $1$ and ${\Bbb H}^2$ is the hyperbolic plane of curvature $-1$.

The systematic study of minimal surfaces in ${\Bbb S}^2\times{\Bbb R}$ and  ${\Bbb H}^2\times{\Bbb R}$ was initiated
 by H. Rosenberg and W. Meeks \cite{mrcmh, Rosenberg} and has been very active since then.
The existence of an associate family for simply connected minimal surfaces in ${\Bbb S}^2\times{\Bbb R}$ and  ${\Bbb H}^2\times{\Bbb R}$ was proved in \cite{Daniel2}.

It has been shown in Bayrakdar et al. \cite{Bayrakdar1}, the Gaussian curvature of a
surface corresponding to a first-order ODE is given by certain Burgers' equations
and it is possible to obtain two-dimensional spaces of constant curvature
from some integrable PDEs.

The main idea for writing present paper is raised from the papers \cite{Bayrakdar3}, written by  T. Bayrakdar and A. A. Ergin.
They proved that a surface corresponding to a first-order ODE is
minimal in three-dimensional Riemannian manifold which is determined
by the first prolongation of $y' = p(x, y)$, if and only if $p_{yy} = 0$. 
Hence any linear first-order ODE characterizes a minimal surface which
is not necessarily totally geodesic.
Z. O. Bayrakdar and T. Bayrakdar used the same idea and got  good results \cite{Bayrakdar2}.

In present paper, we prove that any surface corresponding to  linear second-order ODEs $y''=\alpha(x) y +\beta(x)$
where $\alpha$ and $\beta$ are two smooth functions in term of $x$,
as a submanifold is minimal in the class of third-order ODEs corresponding to
the  third-order equation $y'''=\alpha(x)y' + \alpha'(x) y  +\beta'(x)$.
Furthermore, we will show the linear second-order ODE $y''=\pm y+\beta(x)$ is the only case that is defined a minimal surface
and is also totally geodesic.
%
%%%%%%%%%%%%%%%%%%%%%%%%%%%%%%%%%%%%%%%%%%%%%%%%%%%
\section{Geometry of submanifolds via moving frame method}
In this section we consider the method of moving frames to investigation of geometry of submanifolds.
For more details we refer  to \cite{Willmore}.

Assume $N$ be a $n$-dimensional Riemannian manifold equipped with metric $g$. Let $M$ is a $m$-dimensional submanifold of $N$,  locally
imbedded in $N$. The submanifold of the orthogonal frame bundle over $N$,  denoted by $F(N, M)$,
 which includes of adapted frames $\{{\bf e}_A\}, ~ A =1, \ldots, n$ of which $\{{\bf e}_i\},~ i =1, \ldots, m$ are tangent to
$M$ and $\{{\bf e}_\alpha\},~\alpha =m + 1, \ldots, n$ are normal to $M$.  We have the matrix equation
$\widetilde{\bf E}={\bf E}K$ between two adapted frames $\widetilde{\bf E}$ and ${\bf E}$
where $K$ is the matrix
$$
K=\left[
\begin{array}{cccc}
  A & {\rm O}\\[2mm]
    {\rm O} & B
  \end{array}
\right],
$$
where $A\in {\rm O}(m)$ and $B\in{\rm O}(n - m)$.

The matrix of 1-forms $(\omega_B^A)$ is the Levi-Civita connection for $N$  where we have
\begin{eqnarray}
d\omega^A+\omega_B^A\wedge \omega^B&=&0,\\[2mm]
\omega_B^A+\omega_A^B&=&0.
\end{eqnarray}
We can find a connection over $M$ by a  restriction  $\omega_B^A$ to $\omega_j^i$ over $M$. Namely, when we restrict the equation
\begin{eqnarray}
d\omega^A=-\sum_B\omega_B^A\wedge \omega^B,
\end{eqnarray}
to $M$ gives
\begin{eqnarray}\label{fse}
d\omega^i=-\sum_B\omega_B^i\wedge \omega^B=-\sum_j\omega_j^i\wedge \omega^j,
\end{eqnarray}
because $\omega^\alpha=0$ for tangent vectors on $M$. Furthermore $\omega_j^i=-\omega_i^j$, so the matrix $(\omega_j^i)$ is the uniquely
defined Levi-Civita connection corresponding to the metric on $M$
induced from that on $N$.
Takeing the exterior derivative of $\omega^\alpha=0$, leads to
\begin{eqnarray}\label{sse}
0=d\omega^\alpha=-\sum_A\omega_A^\alpha\wedge \omega^A=-\sum_i\omega_i^\alpha\wedge \omega^i.
\end{eqnarray}
From Cartan's lemma it follows that
\begin{eqnarray}
\omega^\alpha_i=\sum_j h_{ij}^\alpha\; \omega^i,
\end{eqnarray}
for smooth functions $h_{ij}^\alpha=h_{ji}^\alpha$.
%
%%%%%%%%%%%%%%%%%%
\section{Submanifolds in  Riemannian manifold}
 Let $(S, g)$ be a  Riemannian
manifold and $U\subset S$ be an open subset. For each arbitrary $p\in U$, we can define an orthonormal
frame  ${\bf E}=({\bf e}_1, {\bf e}_2, {\bf e}_3, {\bf e}_4)$ where ${\bf e}_i\in T_pU$ and its corresponding dual, the coframe
  $\omega = \{\omega^1, \omega^2, \omega^3, \omega^4\}$ where $\omega^i\in T_p^*U$ and we have  $\omega^i({\bf e}_j)=\delta_{ij}$.
 Therefore the metric tensor is defined by
\begin{eqnarray}
g =\sum_{i,j=1}^4\delta_{ij}~ \omega^i\otimes \omega^j.
\end{eqnarray}
By applying the exterior derivatives on these 1-forms we obtain the first structural equations
\begin{eqnarray}\label{first}
d\omega^i=-\sum_{j=1}^4~\theta_j^i\wedge \omega^j, \;\;\;\;\;\;\; 1\leq i \leq4,
\end{eqnarray}
 where the skew-symmetric matrix of 1-forms $\theta=(\theta_j^i)$  is called the $\frak{o}(3, {\Bbb R})$-valued torsion-free
connection.
Now computing the exterior derivative of $\theta_j^i$ gets the second structural equations as following
\begin{eqnarray}\label{first0}
\Omega^i_j=d\theta_j^i+\sum_{k=1}^4~\theta_{k}^i\wedge \theta_j^k,
\end{eqnarray}
 where the skew-symmetric matrix  $\Omega=(\Omega_j^i)$ is called the Riemannian curvature tensor.
 We can rewrite the $\Omega^i_j$ with respect to  coframe $\omega$ as follows
\begin{eqnarray}\label{first00}
\Omega^i_j=\sum_{k<l}R_{jkl}^i~\omega^k\wedge \omega^l.
\end{eqnarray}
Let $(\widetilde{S}, \widetilde{g})$ be an isometrically immersed submanifold
of the surface of $(S, g)$ by inclusion map $\sigma: \widetilde{S}\longrightarrow S$
satisfies $\sigma^*g = \widetilde{g}$ and  $\widetilde{\omega}^4 = 0$.
Suppose that $\widetilde{\bf E} = (\widetilde{\bf e}_1, \widetilde{\bf e}_2, \widetilde{\bf e}_2, {\bf n})$ is an adapted frame
on $S$ with corresponding coframe $\widetilde{\omega} = (\widetilde{\omega}^1, \widetilde{\omega}^2, \widetilde{\omega}^3, \widetilde{\omega}^4)$, such that
 $(\widetilde{\bf e}_1, \widetilde{\bf e}_2, \widetilde{\bf e}_3)$ is tangent
to $\widetilde{S}$ and ${\bf n}$ is normal to $\widetilde{S}$ and ${\bf E}$ and $\widetilde{\bf E}$ have the same orientation.
As a result, the surface $\widetilde{S}$ is equipped with metric tensor
\begin{eqnarray}\label{m}
\widetilde{g}=\sum_{i=1}^4\widetilde{\omega}^i\otimes \widetilde{\omega}^i.
\end{eqnarray}
The orthogonal frames $\widetilde{\bf E}$ and ${\bf E}$ are related by the equation $ \widetilde{\bf E}={\bf E}A$,
 where $A\in SO(4, {\Bbb R})$,  thus we can write the associated connection 1-forms as follows
\begin{eqnarray}\label{cf1}
\widetilde{\theta}=A^{-1} dA + A^{-1} \theta A.
\end{eqnarray}
 Using \eqref{fse}, we can deduce the first structure equations for coframe $\widetilde{\omega}$ by following formula
\begin{eqnarray}\label{cf2}
d\widetilde{\omega}^i=-\sum_{j=1}^4~\widetilde{\theta}_j^i\wedge \widetilde{\omega}^j.\;\;\;\;\;\;\; 1\leq i\leq4
\end{eqnarray}
Take a look at the formula \eqref{sse}, we derive the  second structure equations
\begin{eqnarray}\label{cf4}
d\widetilde{\omega}^i&=&-\sum_{j}~\widetilde{\theta}_j^i\wedge \widetilde{\omega}^j,\;\;\;\;\;\;\; i=1, 2, 3\\ \label{cf5}
0&=&\widetilde{\theta}_1^4\wedge \widetilde{\omega}^1+\widetilde{\theta}_2^4\wedge \widetilde{\omega}^2+\widetilde{\theta}_3^4\wedge \widetilde{\omega}^3,
\end{eqnarray}
on the $\widetilde{S}$. It is necessary to mention the \eqref{cf5} concludes by $\widetilde{\omega}^4=0$.
 An immediate consequences of the equation \eqref{cf5} are
\begin{gather}\label{cf6}
\begin{split}
\widetilde{\theta}_1^4&=h_{11} \widetilde{\omega}^1+h_{12} \widetilde{\omega}^2+h_{13} \widetilde{\omega}^3,\\
\widetilde{\theta}_2^4&=h_{21} \widetilde{\omega}^1+h_{22} \widetilde{\omega}^2+h_{23} \widetilde{\omega}^3,\\
\widetilde{\theta}_3^4&=h_{31} \widetilde{\omega}^1+h_{32} \widetilde{\omega}^2+h_{33} \widetilde{\omega}^3,
\end{split}
\end{gather}
where $h_{ij} = h_{ji}$. These functions are the coefficients of the second fundamental
form
\begin{eqnarray}\label{fund}
{\rm II}=\widetilde{\theta}_1^4 \otimes \widetilde{\omega}^1 + \widetilde{\theta}_2^4 \otimes \widetilde{\omega}^2+ \widetilde{\theta}_3^4 \otimes \widetilde{\omega}^3,
\end{eqnarray}
 Now one can evaluate the  $\widetilde{\Omega}^i_j$ by following formula
\begin{eqnarray}\label{first0000}
\widetilde{\Omega}^i_j=\sum_{k<l}\widetilde{R}_{jkl}^i~\widetilde{\omega}^k\wedge \widetilde{\omega}^l,
\end{eqnarray}
for metric $\widetilde{g}$ on $\widetilde{S}$.
According to Gauss formula we have
\begin{eqnarray}\label{cf7}
R_{ijrs}-\widetilde{R}_{ijrs}=h_{is}h_{jr}-h_{ir}h_{js}.
\end{eqnarray}
Also in an orthonormal frame we have
\begin{eqnarray}\label{cf8}
\sum_{k}\delta_{ik}R_{jrs}^k=R_{ijrs}.
\end{eqnarray}
We shall define the Weingarten operator ${\cal A}: TM\longrightarrow TM$,  for each $X, Y\in  T_pM$, by
\begin{eqnarray}\label{cf9}
\langle {\cal A}X, Y\rangle =\langle {\rm II}(X, Y), {\bf n}\rangle _N,
\end{eqnarray}
componentwise this means that
\begin{eqnarray}\label{cf10}
{\cal A}=\sum_{i, j} h_{ij}\widetilde{\omega}^i\otimes \widetilde{\bf e}_j.
\end{eqnarray}
In fact, often we shall not distinguish between ${\cal A}$ and the second fundamental
tensor in the direction of ${\bf n}$, that is, the map $\langle {\rm II}(, ), ~{\bf n}\rangle _N : TM\times TM\longrightarrow {\Bbb R}$.
The $k$-th mean curvatures of the hypersurface in the direction
of ${\bf n}$  are given by
\begin{eqnarray}\label{cf11}
H_k=\binom{m}{k}^{-1} S_k,
\end{eqnarray}
where $S_0=1$ and, for $1\leq k\leq m$, $S_k$ is the $k$-th elementary symmetric function of
the eigenvalues of ${\cal A}=(h_{ij})$. In particular $H_1=H$ is the mean curvature that
\begin{eqnarray}\label{meancu2}
H=\frac{1}{3}{\rm tr}~{\cal A},
\end{eqnarray}
 and $H_m$ is the Gauss-Kronecker curvature that equals to
 \begin{eqnarray}\label{meancu1}
H_m={\rm det}~{\cal A},
\end{eqnarray}
and $H_2$ is strictly related to the scalar curvature of $M$, \cite{Luis}.

We remind the reader that the surface $\widetilde{S}$ is said to be {\it totally geodesic} if the second fundamental form
identically vanishes on $\widetilde{S}$ and is said to be {\it minimal} if $H = 0$.
%
%%%%%%%%%%%%%%%%%%%%%%%%%%%%%%%%%%%%%%%%%%%%%%%
\section{The class of third-order ODEs as a Riemannian geometry}
Geometrically, one can consider the third-order ODE with the following form
\begin{eqnarray}\label{3ordeq}
\dfrac{d^3y}{d x^3}=f(x, y, y', y''),
\end{eqnarray}
 as a submanifold ${\cal S}$ in the third-order jet bundle
${\rm J}^3$, which has local coordinates
$$\Upsilon=\{(x, y, p, q, r)\in{\rm J}^3: p=y', q=y'', r=y'''\},$$
and it is denoted by the zero set of the function $F(x, y, p, q, r)=r-f(x, y, p, q)$, that means
${\cal S}=F^{-1}(0)$ is presented on ${\rm J}^2$ as the graph of the function $r=f(x, y, p, q)$.

From the geometric theory of such equations, it follows that there exists
a collection of independent $1$-forms as a coframe
\begin{eqnarray} \nonumber
\omega^1&=&dx,\\ \nonumber
\omega^2&=&dy-p\; dx,\\ \label{cofram}
\omega^3&=&dp-q\; dx,\\ \nonumber
\omega^4&=&dq-f\; dx,\nonumber
\end{eqnarray}
on ${\Bbb R}^5$ with coordinates $(x, y, p, q, r)$.
The third prolongation of a solution curve of differential equation \eqref{3ordeq}
is a curve on ${\cal S}$  represented by 3-jet of a smooth section
$\sigma(x, y, p)=(x, y, p, q, r)$ of the trivial bundle $\pi:{\Bbb R}\times {\Bbb R}\longrightarrow{\Bbb R}$ on which
the contact forms $\omega^2$, $\omega^3$ and $\omega^4$ vanish.
Since the local coframe $\omega=\{\omega^1, \omega^2, \omega^3, \omega^4\}$ is dual to the frame
of the vector fields $\Omega=\{{\bf e}_1, {\bf e}_2, {\bf e}_3, {\bf e}_4\}$ that means $\omega^i({\bf e}_j)=\delta_{ij}$
therefore we have
\begin{eqnarray} \label{fram}
{\bf e}_1=\partial_x+p\partial_y+q\partial_p+f\partial_q,\;\;\;\;\;\;\;\;{\bf e}_2=\partial_y,\;\;\;\;\;\;\;\;\;\;{\bf e}_3=\partial_p,\;\;\;\;\;\;\;\;{\bf e}_4=\partial_q.
\end{eqnarray}
The Riemannian metric
\begin{eqnarray}
g =\sum_{i,j=1}^4\delta_{ij}~ \omega^i\otimes \omega^j,
\end{eqnarray}
on ${\cal S}$ is given in coordinates $(x, y, p, q)$ as
\begin{eqnarray}\label{metr}
ds^2 = (1+p^2 + q^2 + f^2)dx^2 - 2pdxdy  -2q dx dp  - 2fdxdq + dy^2 +dp^2 +dq^2.
\end{eqnarray}
Differentiating the coframe \eqref{cofram}, we have
\begin{gather}\label{dcofram}
\begin{split}
d\omega^1&= 0,\\
d\omega^2&= \omega^1 \wedge \omega^3,\\
d\omega^3&= \omega^1 \wedge \omega^4, \\
d\omega^4&= f_y ~\omega^1 \wedge \omega^2 + f_p ~\omega^1 \wedge \omega^3 + f_q~\omega^1 \wedge \omega^4.
\end{split}
\end{gather}
Now using the \eqref{first} formula, we can write
$$
\left[
\begin{array}{cccc}
   d\omega^1 \\[2mm]
    d\omega^2 \\[2mm]
   d\omega^3 \\[2mm]
  d\omega^4  \\[2mm]
  \end{array}
\right]
=
-\left[
  \begin{array}{cccc}
   0 & \theta^1_2 &\theta^1_3 & \theta^1_4\\[2mm]
    \theta_1^2 &  0& \theta^2_3 & \theta^2_4 \\[2mm]
    \theta_1^3 & \theta_2^3 & 0 & \theta^3_4 \\[2mm]
    \theta_1^4 & \theta_2^4 & \theta_3^4 & 0 \\[2mm]
  \end{array}
\right]
~\wedge~
\left[
\begin{array}{cccc}
   \omega^1 \\[2mm]
    \omega^2 \\[2mm]
   \omega^3 \\[2mm]
  \omega^4  \\[2mm]
  \end{array}
\right],
$$
where $\Theta=(\theta_j^i)$ is a ${\frak o}(4, {\Bbb R})$-valued torsion-free connection is  a antisymmetric matrix, where
\begin{gather}\label{arrays}
\begin{split}
\theta_2^1&=-\theta_1^2= -\frac{1}{2}\ \omega^3 - \frac{1}{2}f_y\ \omega^4, \\
\theta_3^1&=-\theta_1^3= -\frac{1}{2}\ \omega^2 -\frac{1}{2}(1+f_p)\ \omega^4,\\
\theta_4^1&=-\theta_1^4=- \frac{1}{2}f_y\ \omega^2-\frac{1}{2}(1+f_p)\ \omega^3-f_q\ \omega^4, \\
\theta_3^2&=-\theta_2^3=-\frac{1}{2}\ \omega^1, \\
\theta_4^2&=-\theta_2^4= \frac{1}{2}f_y\ \omega^1,\\
\theta_4^3&=-\theta_3^4=-\frac{1}{2}(1- f_p)\ \omega^1.
\end{split}
\end{gather}
Now according to \eqref{first00}, the components of the curvature $2$-form are
\begin{gather*}
\begin{split}
\Omega_2^1&= -\frac{1}{2}\Big[\frac{1}{2}(-1+3f^2_y)\ \omega^1 \wedge \omega^2 + \frac{1}{2}f_y(1+3f_p)\ \omega^1 \wedge \omega^3 + f_{yy}\ \omega^2 \wedge \omega^4\\ &~~~~~~+f_{yp}\ \omega^3 \wedge \omega^4 +\Big(\frac{1}{2}+2f_yf_q-\frac{f_p}{2}
 +pf_{yy}+f_{xy}+ff_{yq}+qf_{yp}\Big)\ \omega^1 \wedge \omega^4\Big],\\
\Omega_3^1&= \frac{1}{2}\Big[- \frac{1}{2}f_y(1+3f_p)\ \omega^1 \wedge \omega^2 - (1+f_p+\frac{3}{2}f_p^2)\ \omega^1 \wedge \omega^3\\
&~~~~~~-(2f_pf_q+pf_{yp}+ff_{pq}+\frac{1}{2}f_y+f_{xp}) \ \omega^1 \wedge \omega^4 -f_{yp}\ \omega^2 \wedge \omega^4 - f_{pp}\ \omega^3 \wedge \omega^4\Big],  \\
\Omega_4^1&= -\frac{1}{2} \Big(2f_yf_q+\frac{1}{2}(1-f_p)+pf_{yy}+f_{xy}+ff_{yq}+qf_{yp}\Big)\ \omega^1 \wedge \omega^2 -
 \frac{1}{2}\Big(2f_pf_q+pf_{yp}+ff_{pq}\\
&~~~~~~ +qf_{pp}+\frac{1}{2}f_y+f_{xp}\Big)\ \omega^1 \wedge \omega^3  -\frac{1}{2}\Big(\frac{3}{2}-\frac{1}{2}f_p^2-\frac{1}{2}f_y^2+2f_q^2+f_p+2qf_{pq}+2ff_{qq}
+2pf_{yq}
 \\
 &~~~~~~+2f_{xq}\Big)\ \omega^1 \wedge \omega^4-\frac{1}{2}f_{yq}\ \omega^2 \wedge \omega^4 -\frac{1}{2}f_{pq}\ \omega^3 \wedge \omega^4,\\
\Omega_3^2 &=\frac{1}{2}\Big[ \frac{1}{2}\ \omega^2 \wedge \omega^3 + \frac{1}{2}f_y\ \omega^2 \wedge \omega^4-\frac{1}{2}(1+f_p)\ \omega^3 \wedge \omega^4\Big],\\[2mm]
\Omega_4^2 &= \frac{1}{2}\Big[-f_{yy}\;\omega^1 \wedge \omega^2 - f_{yp}\ \omega^1 \wedge \omega^3 - f_{yq}\ \omega^1 \wedge \omega^4 +
 \frac{1}{2}f_y\;\omega^2 \wedge \omega^3 + \frac{1}{2}f_y^2\;\omega^2 \wedge \omega^4\\
&~~~~~~-\Big( f_q-\frac{1}{2}f_y-\frac{1}{2}f_yf_q\Big)\;\omega^3 \wedge \omega^4\Big],\\
\Omega_4^3 &= \frac{1}{2}\Big[ -f_{yp}\ \omega^1 \wedge \omega^2 - f_{pp}\ \omega^1 \wedge \omega^3 - f_{pq}\ \omega^1 \wedge \omega^4 -\frac{1}{2}(1+f_p)\ \omega^2 \wedge \omega^3+\frac{1}{2}\Big(f_y(1+f_p)\\
&~~~~~~-2f_q \Big)\ \omega^2 \wedge \omega^4+\frac{1}{2}(1+f_p)^2\ \omega^3 \wedge \omega^4\Big].
\end{split}
\end{gather*}
Thus independent components of the Riemann curvature tensor are given below
\begin{eqnarray*}
R^1_{212}&=& \frac{1}{4}\Big(1-3f^2_y\Big),\;\;\;\; R^1_{213}= - \frac{1}{4}f_y(1+3f_p),\\
 R^1_{214}&=& -\frac{1}{2}\Big(\frac{1}{2}+2f_yf_q-\frac{f_p}{2}+pf_{yy}+f_{xy}+ff_{yq}+qf_{yp}\Big),\\
 R^1_{224}&=& - \frac{1}{2}f_{yy},\;\;\;\;\;\;\;\; R^1_{234}= -\frac{1}{2}f_{yp},\;\;\;\\
R^1_{312}&=& - \frac{1}{4}f_y(1+3f_p),\;\;\;\; R^1_{313}=- \frac{1}{2} (1+f_p+\frac{3}{2}f_p^2),\\
 R^1_{314}&=&-\frac{1}{2}(2f_pf_q+pf_{yp}+ff_{pq}+\frac{1}{2}f_y+f_{xp}),\;\;\; R^1_{324}=  -\frac{1}{2}f_{yp},\;\;\;\; R^1_{334}= - \frac{1}{2}f_{pp},\\
R^1_{412}&=& -\frac{1}{2} \Big(2f_yf_q+\frac{1}{2}(1-f_p)+pf_{yy}+f_{xy}+ff_{yq}+qf_{yp}\Big),\\
R^1_{413}&=& \frac{1}{2}\Big(2f_pf_q+pf_{yp}+ff_{pq}+qf_{pp}+\frac{1}{2}f_y+f_{xp}\Big), \\
R^1_{414}&=&  -\frac{1}{2}\Big(\frac{3}{2}-\frac{1}{2}f_p^2-\frac{1}{2}f_y^2+2f_q^2+f_p+2qf_{pq}+2ff_{qq}+2pf_{yq}+2f_{xq}\Big), \\
R^1_{424}&=&-\frac{1}{2}f_{yq}, \;\;\;\;\;\;\; R^1_{434}=-\frac{1}{2}f_{pq}, \\
R^2_{323}&=&\frac{1}{4},\;\;\;\;\;\;\;  R^2_{324}=\frac{1}{4}f_y, \;\;\;\;\;\;\;  R^2_{334}=-\frac{1}{4}(1+f_p),\\
R^2_{412}&=&-\frac{1}{2}f_{yy},\;\;\;\;\;  R^2_{413}=-\frac{1}{2} f_{yp},  \;\;\;\;\;  R^2_{414}=-\frac{1}{2} f_{yq},\;\;\;\;  R^2_{423}= \frac{1}{4}f_y,\;\;\;\;\; R^2_{424}=  \frac{1}{4}f_y^2,\\
R^2_{434}&=&-\frac{1}{4}\Big( f_q-\frac{1}{2}f_y-\frac{1}{2}f_yf_q\Big), \;\;\; R^3_{412}=-\frac{1}{2}f_{yp}, \;\;\; R^3_{413}=-\frac{1}{2} f_{pp}, \;\;\; R^3_{414}=-\frac{1}{2} f_{pq}, \\
\;\;\;\;\; R^3_{423}&=&-\frac{1}{4}(1+f_p),\;\;\;\; R^3_{424}=\frac{1}{4}\Big(f_y(1+f_p)-2f_q \Big), \;\;\;\;\; R^3_{434}=\frac{1}{4}(1+f_p)^2.
\end{eqnarray*}
%%%%%%%%%%%%%%%%%%%%%%%%%%%%%%%%%%%%%%%%%%%%
\section{Geometry of submanifolds}
Let $\widetilde{\cal S}$ is a submanifold  in $\cal S$ determined by the
smooth section
\begin{eqnarray}\label{section}
\sigma : (x, y, p) \longrightarrow (x, y, p, q(x, y, p)),
\end{eqnarray}
corresponding to the equation $\displaystyle{y''=q(x, y, p)}$ such that $dp - qdx$ vanishes
on the first prolongation of an integral curve of \eqref{3ordeq}.
In fact, 3-graph of an integral curve of \eqref{3ordeq} is specified by
$$y'''=q_x + q_yy' + q_pp' = q_x + pq_y + qq_p,$$
lies on $\widetilde{\cal S}$. We can compute the  pullbacks of $\omega^i$ by $\sigma$  as follow
\begin{eqnarray*}
\sigma^* \omega^1&=&\omega^1,\\
\sigma^* \omega^2&=&\omega^2,\\
\sigma^* \omega^3&=&\omega^3,\\
\sigma^* \omega^4&=&  q_y\;\omega^2+q_p\; \omega^3.
\end{eqnarray*}
With a simple calculation we find
\begin{eqnarray*}
\sigma^*(d\omega^4)&=&(q_{xy}+pq_{yy}+qq_{yp}+q_yq_p+q_y\sigma^*(f_q))\;\omega^1\wedge \omega^2+(q_{xp}+q_y+pq_{yp}+qq_{pp}\\
&&~~~+q_p^2+q_p\sigma^*(f_q))\;\omega^1\wedge \omega^3, \\[2mm]
 d(\sigma^*\omega^4)&=&(q_{xy}+pq_{yy}+qq_{yp})\;\omega^1\wedge \omega^2+(q_y+q_{xp}+pq_{yp}+qq_{pp})\;\omega^1\wedge \omega^3+q_p\;\omega^1\wedge \omega^4.
\end{eqnarray*}
The equality $\sigma^*(d\omega^4)=d(\sigma^*\omega^4)$,  leads to $q_p=0$ and $q_y\sigma^*(f_q)=0$. Now we  can consider two different cases:
%%%%%%%%%%
\subsection{Case of $q_y\neq0$}
Since  $\sigma^*(f_q)=0$, that is, $f_q=0$ on the submanifold  $\widetilde{\cal S}$, the induced metric on  $\widetilde{\cal S}$ is readily
found as
\begin{eqnarray}\label{metric*}
\sigma^* ds^2 =\omega^1\otimes\omega^1 +  (1+q_y^2)\;\omega^2\otimes\omega^2+  \;\omega^3\otimes\omega^3.
\end{eqnarray}sur
In  local coordinates $(x, y, p)$,  we can rewrite the metric tensor by $\widetilde{g}$ with following form
\begin{eqnarray*}
\widetilde{g} = [1+p^2(1+q_y)^2+q^2]\;dx^2+  (1+q_y^2)\;dy^2+ \;dp^2-2p(1+q_y^2)\;dxdy -2q\;dxdp.
\end{eqnarray*}
For each point in a coordinate neighborhood $(\widetilde{U}; x, y, p)$ of $T\widetilde{\cal S}$, the tangent space $T\widetilde{\cal S}$  is spanned
by the vector fields $\widehat{{\bf e}}_1=\sigma_*\partial_x, \;\widehat{{\bf e}}_2=\sigma_*\partial_y$ and $\widehat{{\bf e}}_3=\sigma_*\partial_p$,
\begin{eqnarray} \label{fram2}
\widehat{{\bf e}}_1=\partial_x+q_x\partial_q,\;\;\;\;\;\;\;\;\widehat{{\bf e}}_2=\partial_y+q_y\partial_q, \;\;\;\;\;\;\;\;\widehat{{\bf e}}_3=\partial_p.
\end{eqnarray}
Now we can rewrite the  vector fields $\widehat{{\bf e}}_1, \widehat{{\bf e}}_2$ and $\widehat{{\bf e}}_3$
in terms of ${\bf e}_1, {\bf e}_2, {\bf e}_3, {\bf e}_4$ with following form
\begin{eqnarray} \label{fram3}
\widehat{{\bf e}}_1={\bf e}_1-p{\bf e}_2-q{\bf e}_3+(q_x-f){\bf e}_4,\;\;\;\widehat{{\bf e}}_2={\bf e}_2+q_y{\bf e}_4, \;\;\;\widehat{{\bf e}}_3={\bf e}_3.
\end{eqnarray}
To find unit normal vector field to $\widetilde{\cal S}$, we can consider the cross product of
$X, Y , Z\in T{\cal S}$ at a given point of a four-dimensional Riemannian manifold.
The cross product on ${\cal S}$ is defined in terms of the Riemannian metric \eqref{metr}
and the volume form by
\begin{eqnarray}
g(X \times Y \times Z, W) := {\rm vol}_g(X, Y, Z, W), \;\;\; \forall ~W \in T{\cal S},
\end{eqnarray}
where the positive definite matrix $(g_{ij})$ is defined by
$$
(g_{ij})=
\left[
  \begin{array}{cccc}
   1+p^2 + q^2 + f^2 & -p & -q & -f \\
    -p & 1 & 0 & 0 \\
    -q & 0 & 1 & 0 \\
    -f & 0 & 0 & 1 \\
  \end{array}
\right].
$$
Let $\displaystyle X=\sum_{i=1}^4 x_i\; {\bf e}_i, \;\; Y=\sum_{i=1}^4 y_i\; {\bf e}_i$ and
$\displaystyle  Z=\sum_{i=1}^4 z_i\; {\bf e}_i$ in  $T{\cal S}$, using \cite{Mert}  we can define the cross product ,
\begin{gather}\label{crospro}
\begin{split}
X\times Y\times Z&=\left|
  \begin{array}{cccc}
  {\bf e}_1 & {\bf e}_2 & {\bf e}_3 & {\bf e}_4 \\
    x_1 & x_2 & x_3 & x_4 \\
    y_1 & y_2 & y_3 & y_4 \\
    z_1 & z_2 & z_3 & z_4 \\
  \end{array}
\right|\\
&=
\left|
  \begin{array}{cccc}
     x_2 & x_3 & x_4 \\
     y_2 & y_3 & y_4 \\
     z_2 & z_3 & z_4 \\
  \end{array}
\right|\;  {\bf e}_1
+
\left|
  \begin{array}{cccc}
    x_1  & x_3 & x_4 \\
    y_1  & y_3 & y_4 \\
    z_1  & z_3 & z_4 \\
  \end{array}
\right|\; {\bf e}_2
+
\left|
  \begin{array}{cccc}
    x_1 & x_2  & x_4 \\
    y_1 & y_2  & y_4 \\
    z_1 & z_2  & z_4 \\
  \end{array}
\right|\; {\bf e}_3
+
\left|
  \begin{array}{cccc}
    x_1 & x_2 & x_3  \\
    y_1 & y_2 & y_3 \\
    z_1 & z_2 & z_3  \\
  \end{array}
\right|\; {\bf e}_4.
\end{split}
\end{gather}
According to \eqref{crospro}, we can compute
$$\widehat{{\bf e}}_1\times \widehat{{\bf e}}_2\times \widehat{{\bf e}}_3=q_y{\bf e}_2 - {\bf e}_4,$$
and therefore we obtain the unit normal vector field on  $\widetilde{\cal S}$ as follow
\begin{eqnarray}\label{unitnormal}
{\bf n}=\frac{1}{\sqrt{1+q_y^2}}\left(-q_y{\bf e}_2 + {\bf e}_4\right).
\end{eqnarray}
Thus the adapted orthogonal frame $\widetilde{{\bf E}}= (\widetilde{{\bf e}}_1, \widetilde{{\bf e}}_2, \widetilde{{\bf e}}_3, {\bf n})$  on
 $\widetilde{\cal S}$ is
\begin{eqnarray}\label{adaptedframe}
\widetilde{{\bf e}}_1={\bf e}_1, \;\;\;\;\; \widetilde{{\bf e}}_2=\frac{1}{\sqrt{1+q_y^2}} \left({\bf e}_2+q_y{\bf e}_4 \right),\;\;\;\;\;\;\widetilde{{\bf e}}_3={\bf e}_3, \;\;\;\;\;\;\; {\bf n}=\frac{1}{\sqrt{1+q_y^2}}\left(-q_y{\bf e}_2 + {\bf e}_4\right),
\end{eqnarray}
that  frame $\widetilde{{\bf \Omega}}= (\widetilde{\omega}^1, \widetilde{\omega}^2, \widetilde{\omega}^3, \widetilde{\omega}^4)$  on
 $\widetilde{\cal S}$ is

\begin{eqnarray}\label{orthocoframe}
\widetilde{\omega}^1={\omega}^1, \;\;\;\; \widetilde{\omega}^2=\frac{1}{\sqrt{1+q_y^2}} \left({\omega}^2+q_y{\omega}^4 \right),
\;\;\;\;\;\widetilde{\omega}^3={\omega}^3, \;\;\;\;
\widetilde{\omega}^4=\frac{1}{\sqrt{1+q_y^2}}\left(-q_y{\omega}^2 + {\omega}^4\right).
\end{eqnarray}
Now since we have $\sigma^* \omega^4= q_y\;\omega^2$ on the submanifold $\widetilde{\cal S}$, we will have
$\widetilde{\omega}^4=0$ on it.  Therefore  the
Riemannian metric \eqref{metric*} can be rewritten based on the members of coframe $\widetilde{{\bf \Omega}}$ as follows
\begin{eqnarray}\label{metric**}
 \widetilde{g}=\widetilde{\omega}^1\otimes\widetilde{\omega}^1 +  \widetilde{\omega}^2\otimes\widetilde{\omega}^2+  \widetilde{\omega}^3\otimes\widetilde{\omega}^3,
\end{eqnarray}
on the submanifold $\widetilde{\cal S}$, here $\widetilde{\omega}^2$ is equal to $\sqrt{1+q_y^2}\;\omega^2$.

We have the relation $\widetilde{{\bf E}}={\bf E} A$ between the frames
 ${\bf E}=({\bf e}_1, {\bf e}_2, {\bf e}_3, {\bf e}_4)$ and
 $\widetilde{{\bf E}}= (\widetilde{{\bf e}}_1, \widetilde{{\bf e}}_2, \widetilde{{\bf e}}_3, {\bf n})$
where $A\in{\rm SO}(4, {\Bbb R})$, with following representation
\begin{eqnarray}\label{matrix}
A=
\left[
  \begin{array}{cccc}
   1 & 0 &0 & 0 \\
    0 & \cos(\varepsilon) & 0 & -\sin(\varepsilon) \\
    0 & 0 & 1 & 0 \\
    0 & \sin(\varepsilon) & 0 & \cos(\varepsilon) \\
  \end{array}
\right],
\end{eqnarray}
 is defined by putting $q_y = \tan(\varepsilon)$ for the sufficiently small
values of $\varepsilon$.
Using the formula \eqref{cf1}, the connection matrix $\widetilde{\theta}$ is
equal to
$$
\widetilde{\theta}=
\left[
\begin{smallmatrix}
   0 & \cos(\varepsilon)\;\theta^1_2+ \sin(\varepsilon)\;\theta^1_4 & \theta^1_3 & -\sin(\varepsilon)\;\theta^1_2+ \cos(\varepsilon)\;\theta^1_4 \\[3mm]
   - \cos(\varepsilon)\;\theta^1_2- \sin(\varepsilon)\;\theta^1_4  &  0 & \cos(\varepsilon)\;\theta^2_3- \sin(\varepsilon)\;\theta^3_4  & -d\varepsilon+\theta^2_4 \\[3mm]
    -\theta^1_3 & -\cos(\varepsilon)\;\theta^2_3+ \sin(\varepsilon)\;\theta^3_4  & 0 & \sin(\varepsilon)\;\theta^2_3+ \cos(\varepsilon)\;\theta^3_4  \\[3mm]
   \sin(\varepsilon)\;\theta^1_2-\cos(\varepsilon)\;\theta^1_4 & d\varepsilon-\theta^2_4 & -\sin(\varepsilon)\;\theta^2_3-\cos(\varepsilon)\;\theta^3_4 & 0
    \end{smallmatrix}
\right],
$$
where the 1-forms $\theta_i^j$s are defined in \eqref{arrays}. Since $q_y = \tan(\varepsilon)$, thus
\begin{eqnarray}\label{de}
d\varepsilon=\dfrac{(q_{xy}+pq_{yy})\widetilde{\omega}^1+q_{yy}(1+q_y^2)^{-\frac{1}{2}}\;\widetilde{\omega}^2}{1+q_y^2},
\end{eqnarray}
and by putting $\sin(\varepsilon)=q_y(1+q_y^2)^{-\frac{1}{2}}$ and $\cos(\varepsilon)=(1+q_y^2)^{-\frac{1}{2}}$, we can rewrite the arrays of
the matrice $\widetilde{\theta}$ as follow
\begin{gather}\label{arrays1}
\begin{split}
\widetilde{\theta}_2^1&= -\frac{1}{2}\Big[f_y\sin(\varepsilon)\cos(\varepsilon)\;\widetilde{\omega}^2+(\cos(\varepsilon)+(1+f_p)\sin(\varepsilon))\;\widetilde{\omega}^3\Big],\\
\widetilde{\theta}_3^1&= -\frac{1}{2}\Big[(\cos(\varepsilon)+(1+ f_p)\sin(\varepsilon))\;\widetilde{\omega}^2\Big],\\
\widetilde{\theta}_4^1&= -\frac{1}{2}\Big[f_y \cos(2\varepsilon)\;\widetilde{\omega}^2+(\sin(\varepsilon)-(1+f_p)\cos(\varepsilon))\;\widetilde{\omega}^3\Big],\\
\widetilde{\theta}_3^2&=-\frac{1}{2}\Big[\cos(\varepsilon)-(1- f_p)\sin(\varepsilon) \Big]\widetilde{\omega}^1,\\
\widetilde{\theta_4^2}&=-d\varepsilon +\frac{1}{2}f_y \;\widetilde{\omega}^1,\\
\widetilde{\theta_4^3}&=-\frac{1}{2}\Big[\sin(\varepsilon)+(1- f_p)\cos(\varepsilon) \Big] \widetilde{\omega}^1.
\end{split}
\end{gather}
Looking at \eqref{cf6},  the
coefficients of the second fundamental form are determined by
\begin{eqnarray*}
&&h_{11}=0, \\
&& h_{12}=h_{21}=\frac{1-q_y^2}{1+q_y^2}\; (q_x+pq_y)_y, \\
 &&h_{13}=h_{31}=-\frac{1}{2\sqrt{1+q_y^2}},  \\
 &&h_{22}=-\frac{q_{yy}}{(1+q_y^2)^{\frac{3}{2}}},  \\
&&h_{23}=h_{32}=0, \\
&&h_{33}=0,  \\
\end{eqnarray*}
here ${\cal A}=(h_{ij})$.
Now Using the formula \eqref{fund}, the second fundamental form is
\begin{gather}\label{secfun}
\begin{split}
{\rm II}&=-\left(\frac{1-q_y^2}{1+q_y^2}\right)(q_x+pq_y)_y\;\widetilde{\omega}^1\otimes\widetilde{\omega}^2-\left(\frac{q_y-(q_x+pq_y)_p}{\sqrt{1+q_y^2}}\right)\;\widetilde{\omega}^1\otimes\widetilde{\omega}^3  -\dfrac{q_{yy}}{(1+q_y^2)^{3/2}}\;\widetilde{\omega}^2 \otimes\widetilde{\omega}^2
\end{split}
\end{gather}
%%%%%%%%%%%%%%%%%%%
Then using \eqref{meancu1} and \eqref{meancu2}, we have
\begin{eqnarray}
{\cal K}_e&=&\frac{q_{yy}}{4(1+q_y^2)^{\frac{5}{2}}},\\
H&=&-\frac{q_{yy}}{3(1+q_y^2)^{\frac{3}{2}}}.
\end{eqnarray}
Using above findings we can summarize following theorems:
\begin{theorem}
A submanifold $\widetilde{\cal S}\subset \cal S$ which is determined by the section
\eqref{section} is minimal if and only if
\begin{eqnarray}
q_{yy}=0.
\end{eqnarray}
In this case, the Gauss-Kronecker curvature of a minimal surface  equals to zero.
\end{theorem}
%%%%%%%%%%%
\begin{theorem}
Let  $\widetilde{\cal S}$ be a submanifold corresponding to linear second-order differential equation
\begin{eqnarray}
y''=\alpha(x) y +\beta(x)
\end{eqnarray}
where $\alpha$ and $\beta$ are two smooth functions in term of $x$, then $\widetilde{\cal S}$ determines a minimal surface in four-dimensional manifold corresponding to
the following third-order equation
\begin{eqnarray}
y'''=\alpha(x)y' + \alpha'(x) y +\beta'(x),
\end{eqnarray}
\end{theorem}
\begin{theorem}
A submanifold $\widetilde{\cal S}\subset \cal S$, which is defined by the section
$$(x, y, p) \mapsto (x, y, p, \alpha(x) y +\beta(x)),$$
 is totally geodesic if and only if $\alpha(x)=\pm1$.
\end{theorem}
%%%%%%%%%%%%%%%%%%%%%%%%%
%%%%%%%%%%%%%%%%%%%%%%%%%%

\subsection{Case of $q_y=0$}
In this case, we want to consider the submanifold determined with $\sigma^*\omega^4=0$,  on the submanifold  $\widetilde{\cal S}$.
Since $\sigma^* \omega^4=  q_y\;\omega^2$, this case concludes to $q_y=0$.
Thus the matrix \eqref{matrix} reduces to the identity
and then we have $\widetilde{\theta}=\theta$.
According to above discussion $ \omega^4$ equals to zero on the submanifold $\widetilde{\cal S}$.
\begin{gather}\label{arrays02}
\begin{split}
\theta_2^1&=-\theta_1^2= -\frac{1}{2}\ \omega^3,  \\
\theta_3^1&=-\theta_1^3= -\frac{1}{2}\ \omega^2,\\
\theta_4^1&=-\theta_1^4=- \frac{1}{2}f_y\ \omega^2-\frac{1}{2}(1+f_p)\ \omega^3, \\
\theta_3^2&=-\theta_2^3=-\frac{1}{2}\ \omega^1 ,\\
\theta_4^2&=-\theta_2^4= \frac{1}{2}f_y\ \omega^1,\\
\theta_4^3&=-\theta_3^4=-\frac{1}{2}(1- f_p)\ \omega^1,
\end{split}
\end{gather}
Therefore
\begin{gather}\label{arrays100}
\begin{split}
\widetilde{\theta}_2^1&= \frac{1}{2}\Big[f_y\;\widetilde{\omega}^2+(1+f_p)\;\widetilde{\omega}^3\Big],\\
\widetilde{\theta}_3^1&= \frac{1}{2}(1+ f_p)\;\widetilde{\omega}^2,\\
\widetilde{\theta}_4^1&= \frac{1}{2}\Big[f_y\;\widetilde{\omega}^2+ \widetilde{\omega}^3\Big],\\
\widetilde{\theta}_3^2&=-\frac{1}{2}(1- f_p) \widetilde{\omega}^1,\\
\widetilde{\theta_4^2}&=-\frac{1}{2}f_y \;\widetilde{\omega}^1,\\
\widetilde{\theta_4^3}&=-\frac{1}{2}\widetilde{\omega}^1.
\end{split}
\end{gather}
That means $h_{11}=h_{22}=h_{33}=h_{23}=0$ and $h_{12}=\frac{1}{2}f_y,\; h_{13}=\frac{1}{2}$, and therefore $H={\cal K}_e=0$
and the second fundamental form is
\begin{eqnarray}\label{secfun2}
{\rm II}=\frac{1}{2} (q_x+pq_y)_y\;\widetilde{\omega}^1\otimes\widetilde{\omega}^2
+\frac{1}{2}\;\widetilde{\omega}^1\otimes\widetilde{\omega}^3,
\end{eqnarray}
\begin{theorem}
A submanifold $\widetilde{\cal S}$ of  $\cal S$ which is determined by the section
\eqref{section} is a minimal but not totally geodesic.
Its Gauss-Kronecker curvature equals to zero.
\end{theorem}
%%%%%%%%%%%%%%%%%%%%%%%%%%%%%%%%%%%%%%
%\section*{Acknowledgement}
%
%The second  author (R. Bakhshandeh-Chamazkoti) acknowledges the funding support
%of Babol Noshirvani University of Technology under Grant No. BNUT/391024/1400.
%%%%%%%%%%%%%%%%%%%%%%%%%%%%%%%%%%%%%%

\end{document}